\documentclass[a4paper,10pt]{article}
  \usepackage{amsmath,amsfonts}

  \title{On a paper by Y.~G.~Zarhin}
  \date{12 May 2011}
  \author{Elmer Rees}

\begin{document}
 \maketitle

In a recent paper \cite{[Z]} it was proved that each member of a
naturally defined family of  linear maps $ \mathbb{C}^n \to
\mathbb{C}^n $ had co-rank one. \S 1 of this note presents  a direct
proof of Zarhin's result about complex polynomials with distinct
roots; it is rather similar to that of the Appendix to \cite{[Z]} by
V.~S.~Kulikov but I give explicit constants. In \S 2 we discuss the
case of a polynomial with multiple roots.

\bigskip
\noindent {\bf \large \S 1. Zarhin's result} \medskip

First we recall the main result of \cite{[Z]}. Let $f$ be a monic
complex polynomial of degree $n$ with $n$ distinct roots
$\alpha_1,\alpha_2,\ldots,\alpha_n$. Define
$$M(f):=(f'(\alpha_1),f'(\alpha_2),\ldots,f'(\alpha_n)),$$ then
the derivative $dM_f$ of $M$ for each such  $f$ has rank $n-1.$

Since the map ${\bf \alpha} = (\alpha_1,\alpha_2,\ldots,\alpha_n)
\to f(x)=(x-\alpha_1) (x-\alpha_2)\ldots(x-\alpha_n)$ is a regular
($n!$-sheeted) covering at points ${\bf \alpha}$ where the
$\alpha_i$ are distinct, the study of the rank of the derivative of
$M$ can be done equivalently at $f$ or ${\bf \alpha}$. We do so at
${\bf \alpha}$.

The $n \times n$ matrix $T$ of the derivative map is given by
$$T_{ij}=\frac{\partial}{\partial\alpha_j}f'(\alpha_i)=\frac{\partial}{\partial\alpha_j}
(\alpha_i-\alpha_1)(\alpha_i-\alpha_2)\ldots
\widehat{(\alpha_i-\alpha_i)}\ldots (\alpha_i-\alpha_n)$$ \vspace{-5
mm}
$$\begin{array}{cllc}

\hspace{5
mm}=&-(\alpha_i-\alpha_1)(\alpha_i-\alpha_2)\ldots\widehat{(\alpha_i-\alpha_j)}\ldots\widehat{(\alpha_i-\alpha_i)}\ldots
(\alpha_i-\alpha_n)&i\neq j \\
T_{ii}=&\sum_j(\alpha_i-\alpha_1)(\alpha_i-\alpha_2)\ldots\widehat{(\alpha_i-\alpha_j)}\ldots\widehat{(\alpha_i-\alpha_i)}\ldots
(\alpha_i-\alpha_n).&

\end{array}$$
Our proof of Zarhn's result shows that the matrix $T$ has some
remarkable properties and so it might be of independent interest.

We will simplify the notation by writing
$$ f_k(x)= \frac{f(x)}{(x-\alpha_k)} = - \frac{\partial f}{\partial
\alpha_k}$$ and similarly when $k \neq \ell$,$$ f_{k \ell}(x)=
\frac{f(x)}{(x-\alpha_k)(x-\alpha_{\ell})} \hspace{0.5 cm} {\rm
etc.}$$ Then, for $i \neq j,$ $T_{ij}=f_{ij}(\alpha_i)$.

We let $D_k$ denote the determinant of the matrix obtained from $T$
by deleting the $k^{th}$ row and the $k^{th}$ column and $\Delta(g)$
denote the discriminant of a polynomial $g$.

 We note that the sum of the columns of $T$  is
zero and so rank($T$) $<n.$ Since the discriminant of a polynomial
with distinct roots is non-zero, the proof will be completed by

\medskip \noindent {\bf Proposition} For each $k$,
$$D_k=(-1)^{n-1 \choose 2} (n-1)! \Delta(f_k)=(-1)^{n-1 \choose 2} (n-1)!\prod_{1 \leq i
<j<n}(\alpha_i-\alpha_j)^2.$$

 \noindent
{\bf Proof} We prove the result for $k=n$ (that is,  we are
considering the principal minor of $T$) and the proof is, apart from
notation, the same for other values of $k.$ Interchanging both the
$i^{th}$ and $j^{th}$ rows and the $i^{th}$ and $j^{th}$ columns of
$T$ for $1 \leq i <j <n$ interchanges $i$ and $j$ but does not
change the determinant $D_n$ of the principal minor. So $D_n$ is a
symmetric polynomial in $\alpha_1,\alpha_2, \ldots, \alpha_{n-1}.$
If we set $\alpha_i = \alpha_j$ then the $i^{th}$ and $j^{th}$ rows
of $T$ are equal, so if $1 \leq i <j <n$ and $\alpha_i = \alpha_j$
then $D_n=0.$ We recall the well known

\medskip \noindent {\bf Lemma} If $P(x_1,x_2,\dots, x_r)$ is a
symmetric polynomial which vanishes when any pair of the $x$'s are
equal then $p$ is a multiple of $\prod_{1 \leq i <j\leq
r}(x_i-x_j)^2$.

\medskip
So $D_n$ is a multiple of $\prod_{1 \leq i
<j<n}(\alpha_i-\alpha_j)^2$, but they both have total degree
$(n-1)(n-2)$ ($D_n$ because each $T_{ij}$ has total degree $n-2$).
So $D_n=c\Delta(f_n)$  for some constant $c.$

To determine the value of $c,$ we consider each $T_{ij}$ as an
element of the polynomial ring $R[\alpha_1]$ where $R= { \mathbb
C}[\alpha_2,\alpha_3, \ldots, \alpha_n]$. The degrees of the various
$T_{ij}$ as polynomials in $\alpha_1$ are given by

\begin{center}
\begin{tabular}{c|c}

Index  & Degree $T_{ij}$\\
\hline
$i>1,j>1$&1\\
$i>1,j=0$&0\\
$i=1$&$n-2$

\end{tabular}
\end{center}

Moreover, the coefficient of $\alpha_1^{n-2}$ in $T_{1j}$ is $-1$
for $j>1$ and, since the sum of the columns of $T$ is zero, the
coefficient of $\alpha_1^{n-2}$ in $T_{11}$ is $n-1.$ There are no
occurrences of $\alpha_1$ in the first column (except for $T_{11}$)
and so the terms  that contribute $\alpha_1^{2(n-2)}$ to $D_n$ all
come from the product of $T_{11}$ with $D_{1n},$ the determinant of
the matrix $T[1,n]$ obtained by deleting the first and last rows and
columns of $T$. The coefficient of $\alpha_1$ in the entries of
$T[1,n]$ (since these entries are all linear in $\alpha_1$) are
given by
$$\frac{\partial}{\partial \alpha_1}T_{ij} = \frac{\partial}{\partial
\alpha_1}f_{ij}(\alpha_i) =-f_{1ij} \;\; \mbox{for} \;\; i\neq j.$$
But these are precisely the negatives of the off-diagonal entries of
the matrix of type $T$ that one obtains from the polynomial
$f_1(x)=(x-\alpha_2)(x-\alpha_3)\ldots (x-\alpha_n).$ The diagonal
entries also have the same property since the sum of the columns of
$T$ is zero and there are no occurrences of $\alpha_1$ in the
relevant entries of the first column of $T.$ Hence, by induction,
the term in $D_{1n}$ involving $\alpha_1^{n-2} $ is
$\mbox{det}(-I_{n-2})\alpha_1^{n-2}\Delta(f_1)=(-1)^{n-2}\alpha_1^{n-2}\Delta(f_1)$
which proves the Proposition. (It is easy to start the induction
with $n=2.$)

\newpage \noindent {\bf \large \S 2. Multiple roots} \medskip

Now we consider a monic polynomial $f(x) \in {\mathbb C}[x]$ of
degree $n$  with multiple roots. Let $R(f)=\{ \alpha_1,\ldots,
\alpha_r \}$ denote the set of all its roots, $\# R(f)=r$ (where
$r<n$) and $R_k(f)$ the set of roots of $f$ that have multiplicity
exactly $k.$ We order the roots so that their multiplicities are in
decreasing order and suppose that $\#R_1(f)=s;$ clearly $s<n$. The
first $r-s$ rows of the $r \times r$ matrix $M$ are zero, so
$\mbox{rank}(M) \leq s$. Somewhat tentatively, I make the following
conjecture and sketch some of the calculations that support it.

\medskip \noindent
{\bf Conjecture} The rank of $M$ is $s.$
\medskip

Consider an $s \times s$ submatrix $N$ of $M$ formed from a set of
$s$ columns and the last $s$ rows of $M$. We find that if the
determinants $\mbox{det}N$ of $N$ are zero then a pair of roots of
the polynomial are equal. In particular, calculations that I have
carried out suggest that $\mbox{det}N$ is always of the form
$$\pm c \prod (a-b)^t g$$ where the product is over a nonempty set of pairs of
distinct roots $a,b$ of $f$ and $g$ is a polynomial in $\{
\alpha_1,\ldots, \alpha_r \}.$ In various cases, I describe the
powers $t$, the constant $c$  and the polynomials $g$ :

\medskip \noindent
$\bullet\;\;$ Let $f(x)$ have only one root, say, $\alpha_r$ of
multiplicity $1$, then the principal minor $N$ is $ 1 \times1$ and
it is easy to calculate that
$$N=-k_1\frac{f_r(a_r)}{(a_r-a_1)}$$ where $k_1$ is the
multiplicity of the root $\alpha_1$ and (as in \S1) $f_r(x)$ is
$f(x)$ with the factor $(x-a_r)$ omitted.

\medskip \noindent
$\bullet\;\;$ Let $f(x)$ have the root $\alpha_1$ with multiplicity
$k>1$ and  the other roots be $\alpha_2, \ldots,\alpha_r$ all of
multiplicity $1$ then, the principal minor,  $\mbox{det}N$ has
factors $\alpha_1-\alpha_{\ell}$  with index $t=k-1$ and the factors
$\alpha_m-\alpha_{\ell}$ (where $\ell,m >1$) with index $t=2$ and
$c$ is $\pm k(k+1) \ldots (k+r-1).$

\medskip \noindent
$\bullet\;\;$ Let $f(x)$ have the root $\alpha_1$ with multiplicity
$k>1,$ the root $\alpha_2$ with multiplicity $\ell>1$ and the other
roots of multiplicity $1.$ Then, when $r=5,$ the determinant of the
minors have the form indicated but with a non-trivial factor $g.$
The principal minor has $g=\alpha_1+2\alpha_2 -3\alpha_3$ and one of
the other minors has $g=\alpha_1+2\alpha_2 -3\alpha_4$. If both
these $g$ vanish then we have that $\alpha_3 = \alpha_4$ and if some
other factor of the determinant vanishes then, again two of the
$\alpha$ are equal which contradicts our hypothesis.

\medskip
This final calculation seems to indicate that it may be difficult to
verify  the conjecture by a direct calculation.

\medskip \noindent
Department of Mathematics, University of Bristol and \\  School of
Mathematics, University of Edinburgh

    \smallskip \noindent
\texttt{E.Rees@bristol.ac.uk}
\end{document}